\documentclass[11pt,reqno]{amsart}
\usepackage{amsmath}
\usepackage{amssymb}

\newtheorem{thm}{Theorem}[section]
\newtheorem{lem}{Lemma}[section]
\newtheorem{defn}{Definition}[section]

\newtheorem{coro}{Corollary}[section]\numberwithin{equation}{section}
\newtheorem{rmk}{Remark}[section]

\def\pf{{\textit {Proof:} }}
\def\ac{{\textit {Acknowledgement:} }}

\newcommand{\mysection}[1]{\section{#1}\setcounter{equation}{0}}

\newfont{\bb}{msbm10 at 12pt}

\def\R{\hbox{\bb R}}

\newcommand{\bal}{\begin{aligned}}      \newcommand{\eal}{\end{aligned}}
\newcommand{\ba}{\begin{array}}      \newcommand{\ea}{\end{array}}
\newcommand{\bc}{\begin{center}}     \newcommand{\ec}{\end{center}}
\newcommand{\be}{\begin{enumerate}}  \newcommand{\ee}{\end{enumerate}}
\newcommand{\beq}{\begin{eqnarray}}  \newcommand{\eeq}{\end{eqnarray}}
\newcommand{\beQ}{\begin{eqnarray*}} \newcommand{\eeQ}{\end{eqnarray*}}
\newcommand{\bi}{\begin{itemize}}    \newcommand{\ei}{\end{itemize}}
\newcommand{\bt}{\begin{tabular}}    \newcommand{\et}{\end{tabular}}
\newcommand{\bdm}{\begin{displaymath}} \newcommand{\edm}{\end{displaymath}}

\def\qed{\hfill{Q.E.D.}\smallskip}

\newcommand{\ls}{\setlength{\baselineskip}{12pt}
                 \setlength{\parskip}{3mm}}

\title{Hyperbolic positive energy theorem with electromagnetic fields}
\author{Yaohua Wang, Xu Xu}
\address[Yaohua Wang]{Institute of Mathematics, Academy of Mathematics and
System Sciences, Chinese Academy of Sciences, Beijing 100080, PR China
}
\email{wangyaohua@amss.ac.cn}
\address[Xu Xu]{School of Mathematics and Statistics, Wuhan University, Wuhan 430072, PR China
}
\email{xuxu2@whu.edu.cn}

\date{}

\begin{document}
\maketitle

\begin{abstract}
We establish a type of positive energy theorem for asymptotically anti-de Sitter Einstein-Maxwell initial data sets by using Witten's spinoral techniques.
\end{abstract}
\mysection{Introduction}\ls
As a fundamental result in general relativity, the positive mass theorem states that an isolated gravitational system satisfying the dominant energy condition must have nonnegative total mass. Schoen and Yau first gave a rigorous proof of this theorem for asymptotically flat initial data sets \cite{SY1,SY2,SY3}.
Soon later, Witten provided a different proof by using spinoral techniques \cite{Wi}. Parker and Taubes
gave a rigorous proof of the theorem based on Witten's method \cite{PT}.
Later on, Gibbons and Hull successfully generalized Witten's argument
to asymptotically flat initial data sets with electromagnetic fields \cite{GH}.
Gibbons, Hawking, Horowitz and Perry discussed the positive energy theorems for asymptotically flat initial data sets
with black holes \cite{GHHP}.

For certain spacetimes with nonzero cosmological constant, the positive energy theorem should also be true.
Wang \cite{Wax} and  Chru\'{s}ciel and Herzlich \cite{CH} gave the definition of the mass for
asymptotically hyperbolic manifolds in the time symmetric case and the proof of its positivity respectively.
Soon after that, for initial data sets with nonzero second fundamental form, Maerten \cite{Ma} and Chru\'{s}ciel, Maerten and Tod \cite{CMT} established corresponding positive energy theorems, which gave the upper bounds for angular momentum and center of mass. They assumed that the energy-momentum vector $m_{(\mu)}$ was timelike and $m_i$ $ (i=1,2,3)$
may be zero after an asymptotic isometry.
Recently, the first author, Xie and Zhang \cite{WXZ} considered the general case and established corresponding positive energy
theorem for the initial data sets asymptotic to the anti-de Sitter spacetime in polar coordinates
\begin{equation} \label{AdS}
\widetilde{g}_{AdS}=-\cosh^2(\kappa r)dt^2+dr^2+\frac{\sinh^2(\kappa r)}{\kappa^2}(d\theta^2+\sin^2\theta d \psi^2),
\end{equation}
where
$$-\infty < t <\infty,\  0 < r <\infty,\  0\leq \theta < \pi\  ,\  0 \leq \psi < 2\pi.$$
We also refer to \cite{AD} and \cite{ACOTZ} for physical discussions.

This paper  aims at generalizing the results in  \cite{WXZ} to asymptotically anti-de Sitter spacetimes with electromagnetic fields. Similar to \cite{WXZ}, we generalize Witten's spinoral method to a type of asymptotically anti-de Sitter spacetimes with electromagnetic fields. We provide a rigorous proof of the Poincar\'{e} inequality which ensures the existence and uniqueness of the solution for the Dirac-Witten equation.
Similar to \cite{GH}, we give a lower bound of the total energy in terms of the new quantities we define.

This paper is organized as follows: In Section 2, we give the explicit form of imaginary Killing spinors along
time slices in the anti-de Sitter spacetime. These results are contained in \cite{WXZ}. We write them down here just for completeness.
In Section 3, we give the definition of total energy, total momenta,
total electric charge, and total magnetic momenta for asymptotically anti-de Sitter Einstein-Maxwell initial data sets.
In Section 4, we derive the Weitzenb\"{o}ck formula and prove the existence and uniqueness of
the solution for the Dirac-Witten equation.
In Section 5, we prove our positive energy theorem for asymptotically anti-de Sitter Einstein-Maxwell initial data sets and discuss the case with black holes.
A new positive energy theorem for asymptotically anti-de Sitter initial data sets is established also.
Finally, in Section 6, we present our calculations for the Kerr-Newman-AdS spacetime.

\mysection{The anti-de Sitter spacetime}\ls

The anti-de Sitter ($AdS$) spacetime (\ref{AdS}) with negative cosmological constant $\Lambda=-3\kappa^2<0$,
denoted by $(N,\widetilde{g}_{AdS})$, is a static spherically
symmetric solution of the vacuum Einstein equations. The $t$-slice $(\mathbb{H}^3,\breve{g})$ is the hyperbolic
3-space with constant sectional curvature $-\kappa^2$, where $\kappa>0$.

Let the associated orthonormal frame be
 $$\breve{e}_0=\frac{1}{\cosh(\kappa r)}\frac{\partial}{\partial t},\
\breve{e}_1=\frac{\partial}{\partial r}, \ \breve{e}_2=\frac{\kappa}{\sinh(\kappa r)}\frac{\partial}{\partial \theta}, \
\breve{e}_3=\frac{\kappa}{\sinh(\kappa r)\sin\theta}\frac{\partial}{\partial \psi}$$
and denote $\breve{e}^{\alpha}$ as its dual coframe.

For the following application, we fix the following Clifford representation throughout this paper:
\begin{equation*}
\breve{e}_0 \mapsto \begin{pmatrix}\ &\ &1 &\ \\ \ &\ & \ &1\\1&\
&\ &\ \\ \ &1&\ &\ \end{pmatrix}, \ \breve{e}_1 \mapsto \begin{pmatrix}\ &\ &-1 &\ \\ \ &\ & \ &1\\1&\ &\ &\ \\ \ &-1&\ &\
\end{pmatrix},\end{equation*}
\begin{equation}\label{repre}
\breve{e}_2 \mapsto \begin{pmatrix}\ &\ &\ &1 \\ \ &\ & 1 &\ \\\ &-1 &\ &\ \\ -1 &\ &\ &\
\end{pmatrix}, \ \breve{e}_3 \mapsto \begin{pmatrix}\ &\ &\ &\sqrt{-1} \\ \ &\ & -\sqrt{-1} &\ \\\ &-\sqrt{-1} &\ &\ \\ \sqrt{-1} &\ &\ &\
\end{pmatrix}.
\end{equation}

For the anti-de Sitter spacetime, the imaginary Killing spinors $\Phi_0$, defined as
\begin{equation} \nonumber\\
\nabla^{AdS}_X \Phi_0+\frac{\kappa\sqrt{-1}}{2}X\cdot\Phi_0=0, \ \ \ \ \mbox{$\forall X\in T\mathbb{H}^3$,}
\end{equation}
are all of the form
\begin{equation}\label{ik}
\Phi_0=\begin{pmatrix}u^+e^{\frac{\kappa r}{2}}+u^-e^{-\frac{\kappa r}{2}}\\v^+e^{\frac{\kappa r}{2}}+v^-e^{-\frac{\kappa
r}{2}}\\-\sqrt{-1}u^+e^{\frac{\kappa r}{2}}+\sqrt{-1}u^-e^{-\frac{\kappa r}{2}} \\ \sqrt{-1}v^+e^{\frac{\kappa r}{2}}-\sqrt{-1}v^-e^{-\frac{\kappa
r}{2}}\end{pmatrix}
\end{equation}
along the 0-slice, where
\begin{eqnarray*}
u^+ &=& \lambda_1
e^{\frac{\sqrt{-1}}{2}\psi}\sin\frac{\theta}{2}+\lambda_2
e^{\frac{-\sqrt{-1}}{2}\psi}\cos\frac{\theta}{2},\\
u^-&=&\lambda_3
e^{\frac{\sqrt{-1}}{2}\psi}\sin\frac{\theta}{2}+\lambda_4
e^{\frac{-\sqrt{-1}}{2}\psi}\cos\frac{\theta}{2},\\
v^+&=&-\lambda_3
e^{\frac{\sqrt{-1}}{2}\psi}\cos\frac{\theta}{2}+\lambda_4
e^{\frac{-\sqrt{-1}}{2}\psi}\sin\frac{\theta}{2}, \\
v^-&=&-\lambda_1
e^{\frac{\sqrt{-1}}{2}\psi}\cos\frac{\theta}{2}+\lambda_2
e^{\frac{-\sqrt{-1}}{2}\psi}\sin\frac{\theta}{2}.\end{eqnarray*}
Here  $\lambda_1$, $\lambda_2$, $\lambda_3$, and $\lambda_4$ are
arbitrary complex numbers.

\mysection{Definitions}\ls

 Suppose that $(N,\widetilde{g})$ is a Lorentzian manifold with the metric $\widetilde{g}$ of signature $(-1,1,1,1)$ satisfying the Einstein field
 equations
\begin{equation}
\widetilde{Ric}-\frac{\widetilde{R}}{2}\widetilde{g}+\Lambda \widetilde{g}=T,
\end{equation}
where $\widetilde{Ric}$, $\widetilde{R}$
are the Ricci and scalar curvatures of $\widetilde{g}$ respectively, $T$ is the energy-momentum tensor of matter, and $\Lambda$ is the cosmological
constant.  Let $M$ be a 3-dimensional
spacelike hypersurface in $N$ with the induced metric $g$ and $p$ be the second fundamental form of $M$ in $N$.
$E$ and $B$ are two vector fields on $M$, representing the electric and magnetic fields respectively.
The set $(M,g,p,E,B)$ is called an Einstein-Maxwell initial data set.
In this paper, we will focus on the
following type Einstein-Maxwell initial data set.

\begin{defn} An Einstein-Maxwell initial data set $(M, g, p, E, B)$ is asymptotically $AdS$ of order $\tau >\frac{3}{2}$ if:\\
$(1)$ There is a compact set $K$ such that $M_\infty=M\setminus K$ is diffeomorphic to $\mathbb{R} ^3\setminus B_r$, where $B _r$ is the closed ball of radius $r$ with center at the coordinate origin;\\
$(2)$ Under the diffeomorphism,
$g _{ij}=g(\breve{e} _i,\breve{e} _j)=\delta_{ij}+a_{ij}$, $p _{ij}=p(\breve{e} _i,\breve{e} _j)$,  $E=E^i \breve{e} _i$, $B=B^i \breve{e} _i$
satisfy
\begin{equation}\label{decay condition}
\begin{aligned}
 a_{ij}=O(e^{- \tau \kappa r}),& \ \breve{{\nabla}}_k a_{ij}=O(e^{- \tau \kappa r}), \
\breve{\nabla}_l\breve{{\nabla}}_k a_{ij}=O(e^{-\tau \kappa r}),\\
p_{ij}=&O(e^{- \tau \kappa r}), \
\breve{{\nabla}}_k p_{ij}=O(e^{-\tau \kappa r}),\\
|E|_{\breve{g}}=&O(e^{-2 \kappa r}), \
|B|_{\breve{g}}=O(e^{-2 \kappa r}),\\
B^1=&a e^{-2 \kappa r}+O(e^{-3 \kappa r}),
\end{aligned}
\end{equation}
where
$\breve{{\nabla}}$ is the Levi-Civita connection with respect to the hyperbolic metric
$\breve{g}$ and $a$ is independent of $r$;\\
$(3)$ $div E\in L^1(M)$, $div B\in L^1(M)$, and there is a distance function $\rho_z$
such that $T_{00} e^{\kappa \rho_z}$, $T_{0i} e^{\kappa \rho_z}$ $\in L^1(M)$. \end{defn}

\begin{rmk}
\begin{enumerate}
  \item For simplicity, we just assume that the initial data set has only one end. The case of multi-ends is similar.
  \item For the initial data set $(M, g, p, E, B)$ above, one may get another initial data set $(M, g, p, E', B')$
by duality rotation \cite{Wald},
with the same energy-momentum tensor contributed by electromagnetic filed and $|E'|_{\breve{g}}=O(e^{-2 \kappa r})$,
$|B'|_{\breve{g}}=O(e^{-2 \kappa r})$ and $B'^1=O(e^{-3 \kappa r})$ on the end.
Hence we may assume the original initial data set has the property $B^1=O(e^{-3 \kappa r})$.
\end{enumerate}

\end{rmk}

To get our definitions, recall that the AdS spacetime is just the hyperboloid $\{\eta_{\alpha\beta}y^\alpha
 y^\beta=\frac{3}{\Lambda}\}$ of $\R^{3,2}$ with the metric
\beQ
\begin{aligned}
 \eta_{\alpha\beta}dy^\alpha
dy^\beta= -(dy^0)^2+\sum^3_{i=1}(dy^i)^2-(dy^4)^2 .
\end{aligned}
 \eeQ
The ten Killing vectors
\beQ
\begin{aligned}
 U_{\alpha\beta}=y_\alpha
\frac{\partial}{\partial y^\beta}-y_\beta \frac{\partial}{\partial y^\alpha}
 \end{aligned}
 \eeQ
generate rotations for $\R^{3,2}$. By restricting these vectors to the hyperboloid $\{\eta_{\alpha\beta}y^\alpha
 y^\beta=\frac{3}{\Lambda}\}$ with the induced metric, the Killing vectors of $AdS$ spacetime can be derived, denoted as
 $U_{\alpha\beta}$ also. See appendix for the explicit form of $U_{\alpha\beta}$ along the 0-slice.

Denote
$$\mathcal{E}_i=\breve{\nabla}^j g_{ij}-\breve {{\nabla}}_itr_{\breve{g}}(g)-\kappa(a_{1i}-g_{1i}tr_{\breve{g}}(a)),$$
$$\mathcal{P}_{ki}=p_{ki}-g_{ki}tr_{\breve{g}}(p).$$
Then we can define the following quantities for the asymptotically AdS Einstein-Maxwell initial data set.

\begin{defn}
For the asymptotically $AdS$ Einstein-Maxwell initial data set, the total energy is defined as
\begin{equation*}
E_0=\frac{\kappa}{16\pi}\lim_{r\rightarrow \infty}\int_{S_r}\mathcal{E}_1 U_{40}^{(0)}\breve{\omega},
\end{equation*}
the total momenta are defined as
\begin{equation*}
\begin{aligned}
c_{i}=&\frac{\kappa}{16\pi}\lim_{r\rightarrow \infty}\int_{S_r}\mathcal{E}_1 U_{i4}^{(0)}\breve{\omega},\\
c'_{i}=&\frac{\kappa}{8\pi}\sum_{j=2}^{3}\lim_{r\rightarrow \infty}\int_{S_r}\mathcal{P}_{j1}U_{i0}^{(j)} \breve{\omega},\\
 J_{i}=&\frac{\kappa}{8\pi}\sum_{j=2}^{3}\lim_{r\rightarrow \infty}\int_{S_r}\mathcal{P}_{j1}V_{i}^{(j)} \breve{\omega},
\end{aligned}
\end{equation*}
the total electric charge is defined as
\begin{equation*} q=\frac{1}{4\pi}\lim_{r\rightarrow \infty}\int_{S_r}E^1 \breve{\omega},
\end{equation*}
and the total magnetic momentum is defined as
\begin{equation*}
b_\alpha=\frac{1}{4\pi}\lim_{r\rightarrow
\infty}\int_{S_r}B^1 n^\alpha e^{\kappa r}\breve{\omega},
\end{equation*}
where     $ \  \ \  \ \  \ \  \ \  \ \  \ \
\ \  \ \  \  \breve{\omega}= \breve{e}^2\wedge \breve{e}^3,\  \  \  \ U_{\alpha\beta}=U_{\alpha\beta}^{(\gamma)}\breve{e}_{\gamma},$
$$n^0=1, \ n^1=\sin\theta\cos\psi, \ n^2=\sin\theta\sin\psi, \ n^3=\cos \theta.$$

\end{defn}

\begin{rmk}
For Reissner-Nordstr\"{o}m AdS spacetime with the metric
\beQ
\begin{aligned}
\widetilde{g}=-f dt^2
              +f^{-1}d\bar{r}^2+\bar{r}^2d\Omega^2,
              \ \ f=1-\frac{2M}{\bar{r}}+\frac{Q^2}{\bar{r}^2}+\kappa^2 \bar{r}^2
 \end{aligned}
 \eeQ
and
the field strength tensor
\beQ
\begin{aligned}
F=-\frac{Q}{\bar{r}^2}dt\wedge d\bar{r},
 \end{aligned}
 \eeQ
the total electric charge $q$ is just $Q$.

\end{rmk}

\mysection{The Weitzenb\"{o}ck formula and the Dirac-Witten equation}\ls
Let $(N,\widetilde{g})$ be a spacetime with the metric $\widetilde{g}$ of signature $(-1,1,1,1)$ and
$(M,g,p,E,B)$ be an asymptotically $AdS$ Einstein-Maxwell initial data set of $N$. As $M$ is spin, we can choose a spin structure, and thus a spinor
bundle $\mathbb{S}$ over $M$. This is the same as in \cite{PT}. Denote $\widetilde{\nabla}$, $\nabla$ the Levi-Civita connections with respect to
$\widetilde{g}$, $g$ respectively, and their lifts  to the spinor bundle $\mathbb{S}$.
There is a positive definite Hermitian metric $\langle\cdot, \cdot\rangle$ on $\mathbb{S}$, with respect to which $e_i$ is skew-Hermitian and
$e_0$ is Hermitian\cite{PT,Z99,XD}. Furthermore, $\nabla$ is compatible with  $\langle\cdot, \cdot\rangle$, but $\widetilde{\nabla}$ is not.

For a fixed point $x\in M$, we choose a suitable local orthonormal basis $e_0, e_1, e_2, e_3$,
with $\nabla_{e_i}e_j(x)=0$ for $i,j=1,2,3$ and $\widetilde{\nabla}_{e_0}e_j(x)=0$ for $j=1,2,3$, then
$$\widetilde{\nabla}_{e_i}e_0(x)=p_{ij}e_j,\ \ \widetilde{\nabla}_{e_i}e_j(x)=p_{ij}e_0,$$
where $p_{ij}=\widetilde{g}(\widetilde{\nabla}_{e_i}e_0,e_j)$ is the component of the second fundamental form.
Denote $\{e^\alpha\}$  as the dual coframe of $\{e_\alpha\}$. The two connections on the spinor bundle $\mathbb{S}$ are
related by
$\widetilde{\nabla}_i=\nabla_i-\frac{1}{2}p_{ij}e_0\cdot e_j\cdot$.

We define the imaginary Einstein-Maxwell connection as
$$\widehat{\nabla}_i=\nabla_i-\frac{1}{2}p_{ij}e_0\cdot e_j\cdot +\frac{\kappa\sqrt{-1}}{2}e_i\cdot
   -\frac{1}{2}E\cdot e_i\cdot e_0\cdot -\frac{1}{4}\varepsilon_{jkl}B_je_k\cdot e_l\cdot e_i\cdot ,$$
and define the associated Dirac-Witten operator  as
$$\widehat{D}=\sum_{i=1}^{3}e_i\cdot \widehat{\nabla}_i=D-\frac{1}{2}p_{ii}e_0\cdot -\frac{3\kappa}{2}\sqrt{-1}
  -\frac{1}{2}E\cdot e_0\cdot -\frac{1}{4}\varepsilon_{jkl}B_je_k\cdot e_l\cdot ,$$
where $D=\sum_{i=1}^{3}e_i\cdot \nabla_i$. Then the adjoints of this two operators with respect to $\langle\cdot, \cdot\rangle$ are
\begin{equation*}
\begin{aligned}
\widehat{\nabla}_i^*=&-\nabla_i-\frac{1}{2}p_{ij}e_0\cdot e_j\cdot
  +\frac{\kappa\sqrt{-1}}{2}e_i\cdot -\frac{1}{2}e_i\cdot E\cdot e_0\cdot -\frac{1}{4}\varepsilon_{jkl}B_je_i\cdot e_k\cdot e_l\cdot ,\\
\widehat{D}^*=&D-\frac{1}{2}p_{ii}e_0\cdot +\frac{3\kappa}{2}\sqrt{-1}-\frac{1}{2}E\cdot e_0\cdot
          +\frac{1}{4}\varepsilon_{jkl}B_je_k\cdot e_l\cdot .
\end{aligned}
\end{equation*}

Now we derive the Weitzenb\"{o}ck formulas for $\widehat{D}$ and $\widehat{D}^*$.
\begin{thm}
\begin{equation}\label{weitzenbock formula}
\begin{aligned}
\widehat{D}^*\widehat{D}&=\widehat{\nabla}^{*}\widehat{\nabla}+\widehat{\mathcal{R}},\\
\widehat{D}\widehat{D}^*&=\widehat{\nabla}'^{*}\widehat{\nabla}'+\widehat{\mathcal{R}}',
\end{aligned}
\end{equation}
where
\begin{equation*}
\begin{aligned}
\widehat{\nabla}'=&\nabla_i-\frac{1}{2}p_{ij}e_0\cdot e_j\cdot +\frac{\kappa\sqrt{-1}}{2}e_i\cdot
              -\frac{1}{2}E\cdot e_i\cdot e_0\cdot +\frac{1}{4}\varepsilon_{jkl}B_je_k\cdot e_l\cdot e_i\cdot ,\\
\widehat{\mathcal{R}}=&\frac{1}{2}(\mu-\nu_ie_0\cdot e_i\cdot )+div Ee_0\cdot
                  -div Be_1\cdot e_2\cdot e_3\cdot -\frac{\kappa\sqrt{-1}}{2}\varepsilon_{jkl}B_je_k\cdot e_l\cdot ,\\
\widehat{\mathcal{R}}'=&\frac{1}{2}(\mu-\nu'_ie_0\cdot e_i\cdot )+div Ee_0\cdot
                   +div Be_1\cdot e_2\cdot e_3\cdot -\kappa\sqrt{-1}\varepsilon_{jkl}B_je_k\cdot e_l\cdot ,\\
\end{aligned}
\end{equation*}
with
\begin{equation*}
\begin{aligned}
\mu=\frac{1}{2}\Big(R+\big(&\sum p_{ii}\big)^2-\sum_{i,j}p_{ij}^{2}\Big)+3\kappa^2-|E|^2-|B|^2,\\
\nu_i=&\nabla_jp_{ij}-\nabla_ip_{jj}-2\varepsilon_{ijk}B_jE_k,\\
\nu'_i=&\nabla_jp_{ij}-\nabla_ip_{jj}+2\varepsilon_{ijk}B_jE_k.
\end{aligned}
\end{equation*}
\end{thm}

\pf
By straightforward computation, we have
\begin{equation}\label{weitz 1}
\begin{aligned}
\widehat{D}^*\widehat{D}
=&\nabla^*\nabla+\frac{1}{4}\Big(R+\big(\sum p_{ii}\big)^2+9\kappa^2+|E|^2+|B|^2\Big)\\
   &-\frac{1}{2}\nabla_ip_{jj}e_i\cdot e_0\cdot-\frac{1}{2}\nabla_iE_je_i\cdot e_j\cdot e_0\cdot +E_ie_0\cdot \nabla_i\\
   &-E\cdot e_0\cdot D-\frac{1}{4}\varepsilon_{jkl}\nabla_iB_je_i\cdot e_k\cdot e_l\cdot-\varepsilon_{jkl}B_je_k\cdot \nabla_l \\
   &-\frac{3\kappa\sqrt{-1}}{4}\varepsilon_{jkl}B_je_k\cdot e_l\cdot
     -\frac{1}{2}\varepsilon_{jkl}B_jE_ke_l\cdot e_0\cdot .
\end{aligned}
\end{equation}
Note that
\begin{equation}\label{weitz 2}
\begin{aligned}
\widehat{\nabla}_i^*\widehat{\nabla}_i
=&-\nabla_i\nabla_i+\frac{1}{4}\Big(\sum_{i,j}p_{ij}^2+3\kappa^2+3|E|^2+3|B|^2\Big)\\
 &-div E e_0\cdot +\frac{1}{2}\nabla_{i}p_{ij}e_0\cdot e_j\cdot +\frac{1}{2}\varepsilon_{jkl}B_jE_ke_l\cdot e_0\cdot\\
 &-\frac{\kappa\sqrt{-1}}{4}\varepsilon_{jkl}B_je_k\cdot e_l\cdot -E\cdot e_0\cdot D-\frac{1}{2}\nabla_iE_je_i\cdot e_j\cdot e_0\cdot\\
 &+\frac{1}{4}\varepsilon_{jkl}\nabla_i B_je_k\cdot e_l\cdot e_i\cdot
 +E_ie_0\cdot \nabla_i-\varepsilon_{jkl}B_je_k\cdot \nabla_l.
\end{aligned}
\end{equation}
Combining (\ref{weitz 1}) and (\ref{weitz 2}), we have
\begin{equation*}
\begin{aligned}
\widehat{D}^*\widehat{D}
=&\widehat{\nabla}^*\widehat{\nabla}+\frac{1}{4}\Big(R+\big(\sum p_{ii}\big)^2-\sum_{i,j}p_{ij}^2+6\kappa^2-2|E|^2-2|B|^2\Big)\\
 &-\frac{1}{2}(\nabla_jp_{ij}-\nabla_ip_{jj}-2\varepsilon_{ijk}B_jE_k)e_0\cdot e_i\cdot +div E e_0\cdot \\
 &-\frac{\kappa\sqrt{-1}}{2}\varepsilon_{jkl}B_je_k\cdot e_l\cdot
   -\frac{1}{4}\varepsilon_{jkl}\nabla_iB_j(e_i\cdot e_k\cdot e_l\cdot +e_k\cdot e_l\cdot e_i\cdot ).
\end{aligned}
\end{equation*}
Note that
\begin{equation*}
\begin{aligned}
\varepsilon_{jkl}\nabla_iB_j(e_k\cdot e_l\cdot e_i\cdot +e_i\cdot e_k\cdot e_l\cdot )=4div Be_1\cdot e_2\cdot e_3\cdot ,
\end{aligned}
\end{equation*}
thus
\begin{equation*}
\begin{aligned}
\widehat{D}^*\widehat{D}
=&\widehat{\nabla}^*\widehat{\nabla}+\frac{1}{4}\Big(R+\big(\sum p_{ii}\big)^2-\sum_{i,j}p_{ij}^2+6\kappa^2-2|E|^2-2|B|^2\Big)\\
 &-\frac{1}{2}\left(\nabla_jp_{ij}-\nabla_ip_{jj}-2\varepsilon_{ijk}B_jE_k\right)e_0\cdot e_i\cdot \\
 &+div E e_0\cdot -div Be_1\cdot e_2\cdot e_3\cdot
  -\frac{\kappa\sqrt{-1}}{2}\varepsilon_{jkl}B_je_k\cdot e_l\cdot .
\end{aligned}
\end{equation*}
This proves the first formula in $(\ref{weitzenbock formula})$. The second formula is proved similarly.
\qed

\begin{rmk}
If $\kappa=0$, the formulas in $(\ref{weitzenbock formula})$ are reduced to the  Weitzenb\"{o}ck formulas for the asymptotically flat
Einstein-Maxwell initial data set \cite{GH}.  And if $E=B=0$, the formulas in $(\ref{weitzenbock formula})$ are the same as the
formulas for asymptotically AdS spacetimes \cite{XZ,Wax}.
\end{rmk}

The modified dominant energy condition we impose to prove the positive energy theorem is
\begin{equation}\label{DEC}
\begin{aligned}
\frac{1}{2}\mu\geq\max\bigg\{&\sqrt{\frac{1}{4}|\nu|^2+(div E)^2+(div B)^2}+\kappa|B|,\\
                         &\sqrt{\frac{1}{4}|\nu'|^2+(div E)^2+(div B)^2}+4\kappa|B|\bigg\}.
\end{aligned}
\end{equation}

\begin{rmk}
The modified dominant energy condition $(\ref{DEC})$ is the same as the one used to prove the positive energy theorem
for asymptotically flat manifolds with electromagnetic fields if $\kappa=0$, and is reduced to the standard
dominant energy condition $T_{00}\geq \sqrt{\sum_iT_{0i}^2}$ if $E=B=0$.
\end{rmk}

In order to calculate the boundary term of the Weitzenb\"{o}ck formula, we will define a new connection
and see the difference of two connections on the spinor bundle.
Most of the results presented here are due to the work in \cite{AnD,Z1,XZ}.
These are written down here just for completeness.
Recall that $g=\breve{g}+a$ with
$a=O(e^{-\tau\kappa r}), \
\breve{{\nabla}} a=O(e^{-\tau\kappa r}), \
\breve{{\nabla}}\breve{{\nabla}}a=O(e^{-\tau\kappa r})$.
Orthonormalizing $\breve e_i$ yields
\begin{equation*}
e_i=\breve e_i-\frac{1}{2}a_{ik}\breve e_k+o(e^{-\tau\kappa r}).
\end{equation*}
This provides a gauge transformation
$$\mathcal{A}: SO(\breve g) \rightarrow SO(g)$$
$$\breve e_i \mapsto e_i$$ (and in addition $e_0 \mapsto e_0$)
which identifies also the corresponding spinor bundles.

A new connection is introduced by
$\overline{\nabla}=\mathcal{A}\circ\breve {{\nabla}}\circ\mathcal{A}^{-1}$.
Then we have
\begin{lem} \label{expression for energy-momentum}
$(Prop. 3.2$, \cite{Z2}$)$ Let $(M,g,h)$ be a $3$-dimensional asymptotically AdS initial data set. Then
\begin{equation}\label{asymptotical form for energy-momentum}
\sum_{j,\ j\neq i}Re\langle\phi,e_i\cdot
e_j\cdot({\nabla}_j-\overline{\nabla}_j)\phi\rangle=\frac{1}{4}(\breve{\nabla}^j
g_{ij}-\breve{{\nabla}}_i
tr_{\breve g}(g)+o(e^{-\tau\kappa r}))|\phi|^2,
\end{equation}
for all $\phi \in \Gamma (\mathbb{S})$.\end{lem}

We can extend the imaginary Killing spinors $\Phi_0$ in (2.4) on the
end to the whole $M$ smoothly. Then corresponding to $g$ we may get the spinors
$\overline{\Phi}_0=\mathcal{A}\Phi_0$. Let $\widehat{\overline{\nabla}}_X=\overline{\nabla}_X+\frac{\sqrt{-1}}{2}\kappa X\cdot$,
then
\begin{equation}\label{modified ik}
\widehat{\overline\nabla}_j\overline{\Phi}_0
=\frac{\sqrt{-1}}{4}\kappa a_{jk}e_k\cdot\overline{\Phi}_0+o(e^{-\tau\kappa r})\overline{\Phi}_0.
\end{equation}

For any compact set $K\subset M$, denote $H^1(K, \mathbb{S})$ as the completion of
smooth sections of $\mathbb{S}|_{K}$ with respect to the norm
$$||\phi||^2_{H^1(K,\mathbb{S})}=\int_K \left(|\phi|^2+|\widehat{\nabla}\phi|^2\right)dV_g.$$
And let $H^1(M, \mathbb{S})$ be the completion of
compact supported smooth sections $C_{0}^{\infty}(\mathbb{S})$ with respect to the norm
$$||\phi||^2_{H^1(M,\mathbb{S})}=\int_M \left(|\phi|^2+|\widehat{\nabla}\phi|^2\right)dV_g.$$
Then $H^1(K, \mathbb{S})$ and $H^1(M, \mathbb{S})$ are Hilbert spaces.
Furthermore, we have the following Poincar\'{e} inequality.

\begin{lem}\label{poincare inequality}
There is a constant $C>0$ such that
\begin{equation}
\int_M|\phi|^2dV_g\leq C\int_M|\widehat{\nabla}\phi|^2dV_g,
\end{equation}
for all $\phi\in H^1(M, \mathbb{S})$.
\end{lem}

\pf
As $C_{0}^{\infty}(\mathbb{S})$ is dense in $H^1(M, \mathbb{S})$, we just need to prove the inequality for
$\phi\in C_{0}^{\infty}(\mathbb{S})$. The required inequality is then obtained by continuity.

Separate $M$ into two parts $K$ and $M\setminus K$, where $K$ is a compact set and $M\setminus K$ is the end. Then for any $\phi\in  C_{0}^{\infty}(\mathbb{S})$,
\begin{equation}\label{poincare 1}
\int_K |\phi|^2dV_g\leq C\bigg(\int_K |\widehat{\nabla}\phi|^2dV_g +\int_{\partial(M\setminus K)}|\phi|^2d\sigma\bigg),
\end{equation}
where $C$ is a positive constant. The inequality is obtained by contradiction as follows.
Suppose it is not the case, then for any $n\in \mathbb{N}$, there exists $\phi_n\in H^1(K,\mathbb{S})$ with $||\phi_n||^2_{L^2(K,\mathbb{S})}=1$ such that
$$\int_K|\widehat{\nabla}\phi_n|^2dV_g+\int_{\partial(M\setminus K)}|\phi_n|^2d\sigma\leq \frac{1}{n}.$$
Hence $\{\phi_n\}$ is bounded in $H^1(K,\mathbb{S})$, and there is a subsequence, still denoted by $\{\phi_n\}$,
converging to some $\phi_{\infty}\in H^1(K,\mathbb{S})$ weakly in $H^1(K,\mathbb{S})$. Then
$$||\phi_\infty||_{H^1(K,\mathbb{S})}\leq \liminf_{n\rightarrow \infty}||\phi_n||_{H^1(K,\mathbb{S})}=1.$$
By the Rellich theorem, there is a subsequence, still denoted by $\{\phi_n\}$, converging to $\phi_\infty$
strongly in $L^2(K,\mathbb{S})$, thus $||\phi_\infty||_{L^2(K,\mathbb{S})}=1$. Therefore,
$$||\widehat{\nabla}\phi_\infty||_{L^2(K,\mathbb{S})}=0,$$
and
$$\widehat{D}\phi_\infty=\widehat{\nabla}\phi_\infty=0$$
on $K$. Furthermore, $H^1(K,\mathbb{S})$ is compactly embedded into $L^2(\partial(M\setminus K),\mathbb{S})$ by Rellich theorem,
then $\phi_\infty|_{\partial(M\setminus K)}=0$. Extend $\phi_\infty$ to be zero outside $K$ along $\partial(M\setminus K)$,
then $\phi_\infty\in H^1(M,\mathbb{S})$ and $\widehat{D}\phi_\infty=0$. Thus $\phi_\infty\in C^\infty(int K\cup \partial(M\setminus K))$
by the interior regularity of Dirac-type equation \cite{BC}. Then $\phi_\infty\equiv 0$ by the Weitzenb\"{o}ck formula (\ref{weitzenbock formula}) of $\widehat{D}$
and the modified dominant energy condition (\ref{DEC}), which contradicts $||\phi_\infty||_{L^2(K,\mathbb{S})}=1$.

Under the asymptotical condition (\ref{decay condition}), the volume element $dV_g$ on $M_{\infty}=M\setminus K$ is equivalent to $d\mu=e^{2\kappa r}drd\Theta$,
where $d\Theta$ is the unit volume element on $S^2(1)$. We just need to prove corresponding inequality for the new volume element.
\begin{equation*}
\begin{aligned}
\int_{M_\infty}|\phi|^2e^{2\kappa r}drd\Theta
=&\frac{1}{2\kappa}\int_{\mathbb{R}^3\setminus B(0.R_0)}|\phi|^2d(e^{2\kappa r})d\Theta\\
=&-\frac{1}{2\kappa}\int_{\partial B(0,R_0)}|\phi|^2d\sigma
 -\frac{1}{2\kappa}\int_{\mathbb{R}^3\setminus B(0.R_0)}\partial_r|\phi|^2d\mu.\\
\end{aligned}
\end{equation*}
Note that
\begin{equation*}
\begin{aligned}
&\bigg|\frac{1}{2\kappa}\int_{\mathbb{R}^3\setminus B(0,R_0)}\partial_r|\phi|^2d\mu\bigg|\\
\leq &\frac{1}{2\kappa}\int_{\mathbb{R}^3\setminus B(0,R_0)}\left(2|\langle\widehat{\nabla}_{e_1}\phi,\phi\rangle|+|p||\phi|^2+\kappa|\phi|^2+|E||\phi|^2+|B||\phi|^2\right)d\mu\\
\leq&\frac{1}{2}\int_{\mathbb{R}^3\setminus B(0,R_0)}|\phi|^2d\mu
  +\frac{1}{\kappa}\int_{\mathbb{R}^3\setminus B(0,R_0)}|\langle\widehat{\nabla}_{e_1}\phi,\phi\rangle|d\mu\\
  &+\frac{C}{2\kappa}e^{-\tau\kappa R_0}\int_{\mathbb{R}^3\setminus B(0,R_0)}|\phi|^2d\mu\\
\leq &\Big(\frac{1}{2}+\frac{\epsilon}{\kappa}+\frac{C}{2\kappa}e^{-\tau\kappa R_0}\Big)\int_{\mathbb{R}^3\setminus B(0,R_0)}|\phi|^2d\mu
  +\frac{4}{\epsilon\kappa}\int_{\mathbb{R}^3\setminus B(0,R_0)}|\widehat{\nabla}\phi|^2d\mu,
\end{aligned}
\end{equation*}
thus
\begin{equation*}
\begin{aligned}
\int_{\mathbb{R}^3\setminus B(0,R_0)}|\phi|^2d\mu
\leq &-\frac{1}{2\kappa}\int_{\partial B(0,R_0)}|\phi|^2d\sigma+\frac{4}{\epsilon\kappa}\int_{\mathbb{R}^3\setminus B(0,R_0)}|\widehat{\nabla}\phi|^2d\mu\\
      &+\Big(\frac{1}{2}+\frac{\epsilon}{\kappa}+\frac{C}{2\kappa}e^{-\tau\kappa R_0}\Big)\int_{\mathbb{R}^3\setminus B(0,R_0)}|\phi|^2d\mu.
\end{aligned}
\end{equation*}
Choosing $\epsilon$ small enough and $R_0$ large enough gives
\begin{equation}\label{poincare 2}
\int_{\mathbb{R}^3\setminus B(0,R_0)}|\phi|^2d\mu
\leq -C_1\int_{\partial B(0,R_0)}|\phi|^2d\sigma
      +C_2\int_{\mathbb{R}^3\setminus B(0,R_0)}|\widehat{\nabla}\phi|^2,
\end{equation}
for some constants $C_1>0$, $C_2>0$.
Combining (\ref{poincare 1}) and (\ref{poincare 2}), we get the desired Poincar\'{e} inequality.
\qed

Lemma \ref{poincare inequality} has the following useful corollary.
\begin{coro}\label{equivalent norm}
$\int_M|\widehat{\nabla}\phi|^2dV_g$ is an equivalent norm for $H^1(M, \mathbb{S})$.
\end{coro}

To prove the main result, we shall prove the following key existence and uniqueness theorem for the Dirac-Witten equation.
\begin{lem}\label{Dirac-Witten equation}
Suppose $(M,g,p,E,B)$ is an asymptotically AdS Einstein-Maxwell initial data set of order $\tau>\frac{3}{2}$, then
there exists a unique spinor $\Phi_1\in H^1(M,\mathbb{S})$ such that
$$\widehat{D}(\Phi_1+\bar{\Phi}_0)=0,$$
where $\Phi_0$ is the imaginary Killing spinor defined by $(\ref{ik})$.
\end{lem}

\pf
We follow the arguments of \cite{Z2,XZ} to prove the lemma.
set
$$B(\phi,\psi)=\int_M\langle\widehat{D}\phi,\widehat{D}\psi\rangle,$$
then Corollary \ref{equivalent norm} and the modified dominant energy condition (\ref{DEC}) imply that
 $B(\cdot,\cdot)$ is a coercive bounded bilinear form on $H^1(M,\mathbb{S})$.
Set
$$F(\phi)=-\int_M\langle \widehat{D}\overline{\Phi}_0,\widehat{D}\phi\rangle.$$
By the definition of $\Phi_0$, we have
\begin{equation*}
\begin{aligned}
\widehat{\nabla}_i\overline{\Phi}_0
=&(\nabla_i-\overline{\nabla}_i)\overline{\Phi}_0+\widehat{\overline{\nabla}}_i\overline{\Phi}_0
 -\frac{1}{2}p_{ij}e_0\cdot e_j\cdot \overline{\Phi}_0\\
 &-\frac{1}{2}E\cdot e_i\cdot e_0\cdot \overline{\Phi}_0
 -\frac{1}{4}\varepsilon_{jkl}B_je_k\cdot e_l\cdot e_i\cdot \overline{\Phi}_0.
\end{aligned}
\end{equation*}
The expressions (\ref{ik}), (\ref{decay condition}), (\ref{asymptotical form for energy-momentum}) and (\ref{modified ik}) imply
that $\widehat{\nabla}\overline{\Phi}_0\in L^2(M,\mathbb{S})$ and then $\widehat{D}\overline{\Phi}_0\in L^2(M,\mathbb{S})$,
thus $F(\cdot)$ is a bounded linear functional on $H^1(M,\mathbb{S})$.

By the Lax-Milgram theorem, there is a unique $\Phi_1\in H^1(M,\mathbb{S})$ satisfying
$$\widehat{D}^*\widehat{D}\Phi_1=-\widehat{D}^*\widehat{D}\overline{\Phi}_0$$
weakly. Set $\Phi=\Phi_1+\overline{\Phi}_0$, $\Psi=\widehat{D}\Phi$, then $\Psi\in L^2(M,\mathbb{S})$ and $\widehat{D}^*\Psi=0$ weakly.
By the elliptic regularity for Dirac-type equation \cite{BC}, we have $\Psi\in H^1(M,\mathbb{S})$ and
$\widehat{D}^*\Psi=0$
in the classical sense. Then the Weitzenb\"{o}ck formula (\ref{weitzenbock formula}) and the modified dominant energy condition (\ref{DEC}) imply that
$\widehat{\nabla}'\Psi=0$. We thus have $|\partial_i \ln|\Psi|^2|\leq \kappa+|p|+|E|+|B|$ on the complement of the zero set
of $\Psi$ on $M$.
As a consequence, if there exist $x_0\in M$ such that $\Psi(x_0)\neq 0$, then
$$|\Psi(x)|^2\geq |\Psi(x_0)|^2e^{(\kappa+|p|+|E|+|B|)(|x_0|-|x|)}.$$
This implies that $|\Psi(x)|$ is not in $L^2(M,\mathbb{S})$, which is a contradiction. Thus $\Psi\equiv 0$,
and $\widehat{D}(\Phi_1+\overline{\Phi}_0)=0$.
\qed

\mysection{Positive energy theorems}\ls

Let $\phi=\overline{\Phi}_0+\Phi_1$, where $\Phi_1$ is the solution given by Lemma \ref{Dirac-Witten equation}.
By integrating the Weitzenb\"{o}ck formula (\ref{weitzenbock formula}) for $\widehat{D}$, we get
\begin{equation}\label{integrated form of Weitz formula}
\begin{aligned}
&\int_M|\widehat\nabla\phi|^2\ast 1
 +\int_M\langle\phi,\widehat{\mathcal{R}}\phi\rangle\ast1\\
=&\lim_{r\rightarrow
 \infty}Re\int_{S_r}\langle\phi,\sum_{ j\neq i }e_i\cdot
 e_j\cdot\widehat\nabla_j \phi\rangle\ast e^i\\
=&\ \frac{1}{4}\lim_{r\rightarrow
 \infty}\int_{S_r}(\breve{\nabla}^j
 g_{1j}-\breve
 {{\nabla}}_1tr_{\breve{g}}(g))|\Phi_0|^2\breve{\omega}\\
\ &+\frac{1}{4}\lim_{r\rightarrow
 \infty}\int_{S_r} \kappa(a_{k1}-g_{k1}tr_{\breve{g}}(a))\langle\Phi_0,\sqrt{-1}\breve{e}_k\cdot\Phi_0
 \rangle\breve{\omega}\\
\ &-\frac{1}{2}\lim_{r\rightarrow
 \infty}\int_{S_r}(h_{k1}-g_{k1}tr_{\breve{g}}(h))\langle\Phi_0,\breve{e}_0\cdot\breve{e}_k\cdot\Phi_0\rangle
 \breve{\omega}\\
\ &+\lim_{r\rightarrow
 \infty}\int_{S_r}E^1\langle\Phi_0,\breve{e}_0\cdot\Phi_0
 \rangle\breve{\omega}\\
\ &-\lim_{r\rightarrow
 \infty}\int_{S_r}B^1\langle\Phi_0,\breve{e}_1\cdot\breve{e}_2\cdot\breve{e}_3\cdot\Phi_0
 \rangle\breve{\omega}.
\end{aligned}
\end{equation}

By the Clifford representation (\ref{repre}) and the explicit form (\ref{ik}) of $\Phi_0$, the boundary term
is equal to
\beQ
\begin{aligned}8\pi (\bar\lambda_1, \bar\lambda_2, \bar\lambda_3, \bar\lambda_4)Q
(\lambda_1, \lambda_2, \lambda_3, \lambda_4)^t,
\end{aligned}
 \eeQ
in which the matrix
\begin{equation*}
Q=\begin{pmatrix}
E        &     L \\
\bar{L}^t  &    \hat{E}
\end{pmatrix},
\end{equation*}
where
\begin{equation*}
E=\begin{pmatrix}
E_0-c_3       &      c_1-\sqrt{-1}c_2\\
+b_0-b_3       &      +b_1-\sqrt{-1}b_2\\
\ & \ \\
c_1+\sqrt{-1}c_2    & E_0+c_3\\
+b_1+\sqrt{-1}b_2    & +b_0+b_3\\
\end{pmatrix},
\end{equation*}

\begin{equation*}
\hat{E}=\begin{pmatrix}
E_0+c_3       &      -c_1+\sqrt{-1}c_2\\
-b_0-b_3       &      +b_1+\sqrt{-1}b_2\\
\ & \ \\
-c_1-\sqrt{-1}c_2    & E_0-c_3\\
+b_1-\sqrt{-1}b_2    & -b_0+b_3\\
\end{pmatrix},
\end{equation*}

\begin{equation*}
L=\begin{pmatrix}
c'_{3}       &      -c'_{1}+J_{2}\\
+\sqrt{-1}(q-J_{3})      &      +\sqrt{-1}(c_{2}'+J_{1})\\
\ & \ \\
-c'_{1}-J_{2}   & -c'_{3}\\
+\sqrt{-1}(J_{1}-c_{2}')    & +\sqrt{-1}(J_{3}+ q)\\
\end{pmatrix}.
\end{equation*}

The nonnegativity of the Hermitian matrix $Q$ can be seen from the integrated form of the Weitzenb\"{o}ck
formula (\ref{integrated form of Weitz formula})
and the dominant energy condition $(\ref{DEC})$. \\
Denote
\beQ
\begin{aligned}
&{\bf c}=(c_1,c_2,c_3), \ \ {\bf c'}=(c'_1,c'_2,c'_3),\ \
{\bf J}=({J}_1,{J}_2,{J}_3),\ \ {\bf b}=(b_1,b_2,b_3), \\
&|L|^2=2(|{\bf c'}|^2+|{\bf J}|^2+ q^2),\ \ A=b_0^2+|{\bf c}|^2+|{\bf b}|^2+|{\bf c'}|^2+|{\bf J}|^2+q^2.
\end{aligned}
 \eeQ
When we consider the asymptotically AdS initial data set, we can obtain the nonnegative
 Hermitian matrix $Q$ above with $q={\bf b}=0$. For this case, we have the following theorem.
\begin{thm}
Let $(M,g,h)$ be a $3$-dimensional asymptotically AdS initial data set of the spacetime
$(N,\widetilde{g})$ which satisfies the dominant energy condition. Then we have\\
\begin{equation}\label{b}
 E_0\geq max\Big \{
 \Big(\frac{1}{2}|L|^2+2|{\bf c}|^2\Big)^{\frac{1}{2}}-|{\bf c}|,\Big(\frac{1}{2}|L|^2+2| {\bf c'}\times {\bf J} \big|-F^2\Big)_+^\frac{1}{2}-|{\bf c}|\Big\},
\end{equation}
where
\beQ
\begin{aligned}
F^2&=\sqrt{2}|{\bf c}||L|+\sqrt[4]{2}|{\bf c}|^{\frac{1}{2}}|L|^{\frac{3}{2}}.
\end{aligned}
\eeQ
\end{thm}
\pf Since the Hermitian matrix $Q$ is nonnegative, all the principal minors of $Q$ are nonnegative. From the nonnegativity of the first-order and second-order principal minors of the matrix $Q$, one finds $E_0\geq|{\bf c}|$ and $E_0\geq\frac{1}{2}|L|$.

The sum of the third-order principal minors implies that
 \beQ
 \begin{aligned}
0\leq &2E_0\big(E_0^2-|{\bf c}|^2\big)-E_0|L|^2+4|{\bf c}||{\bf c'}\times {\bf J}|\\
\leq& 2E_0\big(E_0^2-|{\bf c}|^2\big)-E_0|L|^2+2|{\bf c}|\big(|{\bf c'}|^2+|{\bf J}|^2\big)\\
\leq& 2E_0\big(E_0^2-|{\bf c}|^2\big)-E_0|L|^2+4|{\bf c}|E_0^2\\
=&2E_0\big[(E_0+|{\bf c}|)^2-2|{\bf c}|^2-\frac{1}{2}|L|^2\big]. \end{aligned}
 \eeQ
If $E_0>0$, then one has
\begin{equation}\label{first lower bound of E_0}
\begin{aligned}
 E_0\geq \Big(\frac{1}{2}|L|^2+2|{\bf c}|^2\Big)^{\frac{1}{2}}-|{\bf c}|.
\end{aligned}
\end{equation}
If $E_0=0$, then this is a trivial inequality as $E_0\geq|{\bf c}|$ and $E_0\geq\frac{1}{2}|L|$.

Direct computation shows that
\beQ
\begin{aligned}
det Q\leq&\big(E_0^2-|{\bf c}|^2-\frac{1}{2}|L|^2\big)^2-4| {\bf c'}\times {\bf J} \big|^2+8 E_0|{\bf c}|| {\bf c'}\times {\bf J} \big|\\
=&\big(E_0^2+|{\bf c}|^2-\frac{1}{2}|L|^2\big)^2+2|{\bf c}|^2|L|^2-\big(2E_0|{\bf c}|-2| {\bf c'}\times {\bf J} \big|\big)^2\\
\leq&\big(E_0^2+|{\bf c}|^2-\frac{1}{2}|L|^2 + \sqrt{2}|{\bf c}||L|\big)^2+ \sqrt{2}|{\bf c}||L|^3-\big(2E_0|{\bf c}|-2| {\bf c'}\times {\bf J} \big|\big)^2.
\end{aligned}
\eeQ
The inequality (\ref{first lower bound of E_0}) implies $E_0^2+|{\bf c}|^2-\frac{1}{2}|L|^2
+\sqrt{2}|{\bf c}||L|\geq0$. One can verify that
\beQ
\begin{aligned}
det Q \leq&\Big(E_0^2+|{\bf c}|^2-\frac{1}{2}|L|^2
+\sqrt{2}|{\bf c}||L|+\sqrt[4]{2}|{\bf c}|^{\frac{1}{2}}|L|^{\frac{3}{2}}\Big)^2\\
&-\Big(2E_0|{\bf c}|-2| {\bf c'}\times {\bf J} \big|\Big)^2.
\end{aligned}
\eeQ
Thus
$E_0^2+|{\bf c}|^2-\frac{1}{2}|L|^2
+\sqrt{2}|{\bf c}||L|+\sqrt[4]{2}|{\bf c}|^{\frac{1}{2}}|L|^{\frac{3}{2}}
\geq-2E_0|{\bf c}|+2| {\bf c'}\times {\bf J} \big|$
and the inequality follows immediately.
 \qed

\begin{rmk} If  ${\bf c}=0$,  the inequality $E_0\geq\Big(\frac{1}{2}|L|^2+2| {\bf c'}\times {\bf J} \big|-F^2\Big)_+^\frac{1}{2}-|{\bf c}|$
is reduced to $(3.1)$ in \cite{CMT}. Here we remove the assumption that $m_{(\mu)}$ is timelike in \cite{CMT}.

\end{rmk}

\begin{thm}\label{main theorem}
Let $(M,g,p,E,B)$ be a $3$-dimensional asymptotically AdS Einstein-Maxwell initial data set of the spacetime $(N,\widetilde{g})$
which satisfies the modified dominant energy condition $(\ref{DEC})$. Then the nonnegativity of the Hermitian matrix $Q$ implies the following inequality:
\beq
E_0\geq max\Big \{\big(b_0^2+\frac{1}{4}|L|^2\big)^\frac{1}{2},\ \ \Big(\frac{1}{2}(|{\bf c}|^2+|{\bf b}|^2)+\frac{1}{8}|L|^2\Big)^\frac{1}{2},
(A+|{\bf b}|^2+|{\bf c}|^2)^\frac{1}{2} \nonumber\\
-|{\bf b}|-|{\bf c}|,\Big[A-4\sqrt{2}\big(\sum_{i}(b_0 c_{i}+q J_i )^2+|{\bf c'}\times {\bf J}|^2\big)^\frac{1}{2}+F_+^\frac{1}{2}\Big]_+^\frac{1}{2} \Big\},\nonumber
 \eeq
where
\beQ
 \begin{aligned}
f_+=max \{f,0\},
\end{aligned}
 \eeQ
\beQ
\begin{aligned}
F= & 32|{\bf c'}\times {\bf J}|^2+4|{\bf c}\times {\bf J}|^2+36\sum_{i}(b_0 c_{i}+q J_i )^2
+4|{\bf b}\cdot {\bf c}|^2+4|{\bf b}\cdot {\bf J}|^2\\
&+4|{\bf b}\cdot {\bf c}'|^2-8\sqrt{2}\big(\sum_{i}(b_0 c_{i}+q J_i )^2+|{\bf c'}\times {\bf J}|^2\big)^\frac{1}{2}A.
\end{aligned}
\eeQ
\end{thm}
\pf
The nonnegativity of the first-order principal minors ensures $E_0\geq0$.
And from the nonnegativity of the second-order principal minors, one finds
\begin{equation} \label{ine1}
E_0\geq \big(b_0^2+\frac{1}{2}(|{\bf c'}|^2+|{\bf J}|^2+q^2)\big)^\frac{1}{2}
\end{equation}
 and
\begin{equation}\label{ine2}
E_0\geq \Big(\frac{1}{2}(|{\bf c}|^2+|{\bf b}|^2)+\frac{1}{4}\big(|{\bf c'}|^2+|{\bf J}|^2+q^2\big)\Big)^\frac{1}{2}.
\end{equation}

The sum of the third-order principal minors is given, up to a positive constant, by
\beQ
 \begin{aligned}
S=
&E_0(E_0^2-A)+2b_0{\bf b}\cdot {\bf c}+2q{\bf b}\cdot {\bf J}+2\varepsilon_{ijk}c_i c_{j}' J_k.
\end{aligned}
 \eeQ
Using the Cauchy inequality, one derives
\beQ
 \begin{aligned}
S\leq
&E_0(E_0^2-A)+2E_0|{\bf b}||{\bf c}|+2|q||{\bf b}||{\bf J}|+2|{\bf c}||{\bf c'}\times {\bf J}|.
\end{aligned}
 \eeQ
As $|q||{\bf J}|\leq\frac{1}{2}(|q|^2+|{\bf J}|^2)\leq E_0^2$
 and
$|{\bf c'}\times {\bf J}|\leq\frac{1}{2}(|{\bf c'}|^2+|{\bf J}|^2)\leq E_0^2$,
one can obtain
\beQ
 \begin{aligned}
E_0\big((E_0+|{\bf c}|+|{\bf b}|)^2-A-|{\bf c}|^2-|{\bf b}|^2\big)\geq0.
\end{aligned}
 \eeQ
If $E_0>0$, one finds
\begin{equation}\label{third lower bound of E_0}
\begin{aligned}
E_0\geq(A+|{\bf b}|^2+|{\bf c}|^2)^\frac{1}{2}-|{\bf b}|-|{\bf c}|.
\end{aligned}
 \end{equation}
When $E_0=0$, the inequality $(\ref{ine1})$, together with the inequality $(\ref{ine2})$, shows that $Q=0$.
In this case the inequality (\ref{third lower bound of E_0}) becomes trivial.\\
Similarly, we have
\begin{equation}\label{middle inequ in thm 5.2}
\begin{aligned}
 &2b_0{\bf b}\cdot {\bf c}+2q{\bf b}\cdot {\bf J}+2\varepsilon_{ijk}c_i c_{j}' J_k\\
 \leq &2(|{\bf b}|^2+|{\bf c}|^2)^\frac{1}{2}\big(\sum_{i}(b_0 c_{i}+q J_i )^2+|{\bf c'}\times {\bf J}|^2\big)^\frac{1}{2}\\
\leq &2\sqrt{2}E_0\big(\sum_{i}(b_0 c_{i}+q J_i )^2+|{\bf c'}\times {\bf J}|^2\big)^\frac{1}{2}.
\end{aligned}
\end{equation}
Therefore,
\beQ
 \begin{aligned}
S\leq
&E_0(E_0^2-A)+2\sqrt{2}E_0\big(\sum_{i}(b_0 c_{i}+q J_i )^2+|{\bf c'}\times {\bf J}|^2\big)^\frac{1}{2}.
\end{aligned}
 \eeQ
This implies
\begin{equation}\label{positivity of part of the 4th lower bound}
\begin{aligned}
E_0^2\geq A-2\sqrt{2}\big(\sum_{i}(b_0 c_{i}+q J_i )^2+|{\bf c'}\times {\bf J}|^2\big)^\frac{1}{2},
\end{aligned}
\end{equation}
if $E_0>0$. When $E_0=0$, this inequality can be derived by $(\ref{ine1})$ and $(\ref{ine2})$.

The determinant of the matrix $Q$ is
\beQ
 \begin{aligned}
det Q=&\big(E_0^2-A\big)^2
+8E_0(b_0{\bf b}\cdot {\bf c}+q{\bf b}\cdot {\bf J}+\varepsilon_{ijk}c_i c_{j}' J_k)\\
&-4|{\bf c} \times {\bf c'}|^2-4|{\bf c} \times {\bf J}|^2-4|{\bf c'}\times {\bf J}|^2-4b_0^2|{\bf b}|^2-4q^2|{\bf b}|^2\\
&-4\sum_{i}(b_0 c_{i}+q J_i )^2
-4|{\bf b}\cdot {\bf c}|^2-4|{\bf b}\cdot {\bf J}|^2-4|{\bf b}\cdot {\bf c}'|^2\\
&-8b_0\varepsilon_{ijk}b_i c_{j}' J_k
-8q\varepsilon_{ijk}b_i c_{j} c'_k\\
\leq&\big(E_0^2-A\big)^2
+8E_0(b_0{\bf b}\cdot {\bf c}+q{\bf b}\cdot {\bf J}+\varepsilon_{ijk}c_i c_{j}' J_k)-4|{\bf c}\times {\bf J}|^2\\
&-4\sum_{i}(b_0 c_{i}+q J_i )^2
-4|{\bf b}\cdot {\bf c}|^2-4|{\bf b}\cdot {\bf J}|^2-4|{\bf b}\cdot {\bf c}'|^2.
\end{aligned}
 \eeQ
By (\ref{middle inequ in thm 5.2}), one obtains
\beQ
 \begin{aligned}
det Q
\leq&\big(E_0^2-A\big)^2
+8\sqrt{2}E_0^2\big(\sum_{i}(b_0 c_{i}+q J_i )^2+|{\bf c'}\times {\bf J}|^2\big)^\frac{1}{2}-4|{\bf c}\times {\bf J}|^2\\
&-4\sum_{i}(b_0 c_{i}+q J_i )^2
-4|{\bf b}\cdot {\bf c}|^2-4|{\bf b}\cdot {\bf J}|^2-4|{\bf b}\cdot {\bf c}'|^2\\
=&(E_0^2-A+4\sqrt{2}\big(\sum_{i}(b_0 c_{i}+q J_i )^2+|{\bf c'}\times {\bf J}|^2\big)^\frac{1}{2}\big)^2-4|{\bf c}\times {\bf J}|^2\\
&+8\sqrt{2}\big(\sum_{i}(b_0 c_{i}+q J_i )^2+|{\bf c'}\times {\bf J}|^2\big)^\frac{1}{2}A-32|{\bf c'}\times {\bf J}|^2\\
&-36\sum_{i}(b_0 c_{i}+q J_i )^2
-4|{\bf b}\cdot {\bf c}|^2-4|{\bf b}\cdot {\bf J}|^2-4|{\bf b}\cdot {\bf c}'|^2.
\end{aligned}
 \eeQ
Using (\ref{positivity of part of the 4th lower bound}), we get
\beQ
\begin{aligned}
E_0^2\geq A-4\sqrt{2}\big(\sum_{i}(b_0 c_{i}+q J_i )^2+|{\bf c'}\times {\bf J}|^2\big)^\frac{1}{2}+F_+^\frac{1}{2}.
\end{aligned}
\eeQ
The inequality follows immediately.
\qed

\begin{rmk}
Suppose $M$ has inner boundary $\Sigma=\cup\Sigma_{i}^{+}\cup\cup\Sigma_{i}^{-}$, where $\Sigma_{i}^{+}$ and $\Sigma_{i}^{+}$
are future and past trapped surfaces, defined as
\begin{equation*}
\begin{aligned}
\Sigma_{i}^{+}=&\{tr(h)-tr(p|_{\Sigma_i})\geq 0\},\\
\Sigma_{i}^{-}=&\{tr(h)+tr(p|_{\Sigma_i})\geq 0\}.
\end{aligned}
\end{equation*}
If we take $e_3$ as the outer unit normal of $\Sigma$ in $M$ and take the boundary conditon
$e_0\cdot e_3\cdot\phi=\pm\phi$ on $\Sigma_i^{\pm}$, then
$$\int_{\Sigma}\langle\phi,e_3\cdot\widehat{D}\phi+\widehat{\nabla}_{e_3}\phi\rangle$$
is non-positive. In such situation, Lemma $\ref{poincare inequality}$ and Lemma $\ref{Dirac-Witten equation}$
are all valid. Similar arguments appear in \cite{XD}. This verifies Theorem $\ref{main theorem}$ for black holes.
\end{rmk}

\begin{rmk}
If $E_0=0$, there are four linearly independent spinors satisfying $\widehat{\nabla}\phi=0$.
The characterization of the manifold $M$ in such case will be addressed elsewhere.
\end{rmk}

\mysection{Kerr-Newman-AdS case}\ls
In this section, we will calculate our definitions for time slices in the
Kerr-Newman-AdS spacetime \cite{C}.

For Kerr-Newman-AdS spacetime,  the metric in the
Boyer-Lindquist coordinates $(\hat{t}, \hat{r}, \hat{\theta}, \hat{\varphi})$ is
 \beQ
 \begin{aligned}
\widetilde{g}=&-\frac{\Delta_{\hat{r}}}{\rho^2}\Big[d\hat{t}-\frac{a \sin^2\hat{\theta}}{\Sigma}d\hat{\varphi} \Big]^2+\frac{\rho^2}{\Delta_{\hat{r}}}d\hat{r}^2 +\frac{\rho^2}{\Delta_{\hat{\theta}}}d\hat{\theta}^2\\
&+\frac{\Delta_{\hat{\theta}}\sin^2\hat{\theta}}{\rho^2}\Big[a d\hat{t}-\frac{\hat{r}^2+a ^2}{\Sigma}d\hat{\varphi} \Big]^2,
 \end{aligned}
 \eeQ
where
 \beQ
 \begin{aligned}
\Delta_{\hat{r}}&=\big(\hat{r}^2+a^2\big)\big(1+\kappa^2\hat{r}^2\big)-2m\hat{r}+
e^2,\\
\Delta_{\hat{\theta}}&=1-\kappa^2a^2\cos^2\hat{\theta},\\
\rho^2&=\hat{r}^2+a^2\cos^2\hat{\theta},\\
\Sigma&=1-\kappa^2a^2. \end{aligned}
 \eeQ
If we take
 \beQ
\begin{aligned}
e^0&=\frac{\sqrt{\Delta_{\hat{r}}}}{\rho}\Big(d\hat{t}-\frac{a \sin^2\hat{\theta}}{\Sigma}d\hat{\varphi} \Big), \ \ \ \ e^1=\frac{\rho}{\sqrt{\Delta_{\hat{r}}}}d\hat{r},\\
e^2&=\frac{\rho}{\sqrt{{\Delta_{\hat{\theta}}}}}d\hat{\theta}, \ \ \ \
e^3=\frac{\sqrt{\Delta_{\hat{\theta}}}\sin \hat{\theta}}{\rho}\Big(a d\hat{t}-\frac{\hat{r}^2+a ^2}{\Sigma}d\hat{\varphi}\Big),
\end{aligned}
 \eeQ
then the field strength tensor is
\beQ
\begin{aligned}
F=&-\frac{1}{\rho^4}e\big(\hat{r}^2-a^2\cos^2\hat{\theta}\big)
 e^0\wedge e^1-\frac{2}{\rho^4}e\hat{r} a \cos\hat{\theta} e^2\wedge e^3.
\end{aligned}
 \eeQ
Similar to the process in \cite{HT}, after the coordinate transformations
\beQ
\begin{aligned}
&t=\hat{t},\ \ \varphi=\hat{\varphi}+\kappa^2 a \hat{t},\ \
\sinh(\kappa r)\cos {\theta}=\kappa \hat{r}  \cos\hat{\theta}, \\
&\sqrt{\Sigma}\sinh(\kappa r) \sin {\theta}=\kappa\sqrt{\hat{r}^2+a^2}\sin\hat{\theta},
\end{aligned}
 \eeQ
the Kerr-Newman-AdS metric can be written as
$$\widetilde{g}=-\cosh^2(\kappa r)dt^2+dr^2+\frac{\sinh^2(\kappa r)}{\kappa^2}(d\theta^2+\sin^2\theta d
\varphi^2)+a_{\mu \nu}d{x}^{\mu}d{x}^{\nu},$$ where the nonzero components $a_{\mu \nu}$ have the following asymptotic
behaviors£º
  \beQ
 \begin{aligned}
a_{t t}&=\frac{2m \kappa}{\sinh(\kappa r)} B^{-5/2}+O(e^{-3 \kappa r}),\\ a_{t \varphi}&=-\frac{2m a
\kappa}{\sinh(\kappa r)}\sin^2\theta B^{-5/2}+O(e^{-3\kappa r}),\\ a_{\varphi \varphi}&=\frac{2m a^2\kappa
}{\sinh(\kappa r)}\sin^4{\theta}B^{-5/2}+O(e^{-3\kappa r}),\\ a_{rr}&=\frac{2m \kappa}{\sinh(\kappa r)^5}\cosh(\kappa
r)^2B^{-3/2}+O(e^{-7\kappa r}),\\ a_{r \theta}&=-\frac{2m \kappa^2 a^2}{\sinh(\kappa r)^4}\cosh(\kappa
r)\sin{\theta}\cos{\theta}B^{-5/2}+O(e^{-6\kappa r}),\\ a_{\theta \theta}&=\frac{2m a^4 \kappa^3}{\sinh(\kappa
r)^3}\sin^2{\theta}\cos^2{\theta}B^{-7/2}+O(e^{-5\kappa r}),\\ B&=1-a^2\kappa^2\sin^2 \theta.
 \end{aligned}
 \eeQ
Simple calculations show that for $t$-slices, the quantities with the order not higher than $e^{-3\kappa r}$ are
  \beQ
 \begin{aligned}
 a_{11}&=16 m\kappa B^{-3/2}e^{-3 \kappa r}+o(e^{-3\kappa r}),\\
 a_{33}&=16 m a^2\kappa^3 B^{-5/2}\sin^2{\theta}e^{-3 \kappa r}+o(e^{-3\kappa r}),\\
 p_{13}&=p_{31}=24 m a\kappa^3 B^{-5/2}\sin{\theta}e^{-3 \kappa r}+o(e^{-3\kappa r}),
 \end{aligned}
 \eeQ
and
 \beQ
 \begin{aligned}
\mathcal{P}_{31}&=24 m a\kappa^3 B^{-5/2}\sin{\theta} e^{-3 \kappa r}+o(e^{-3\kappa r}),\\
\mathcal{E}_1&=2\kappa a_{11}+\partial_r a_{33}.
\end{aligned}
 \eeQ
We also have
\beQ
 \begin{aligned}
E^1&=4 \kappa^2 e B^{-3/2}e^{-2 \kappa r}+o(e^{-2\kappa r}),\\
B^1=&16 \kappa^3 e a B^{-5/2}\cos \theta e^{-3\kappa r}+o(e^{-3\kappa r}).
\end{aligned}
 \eeQ
Finally we get for Kerr-Newman-AdS spacetime
 \beQ
 \begin{aligned}
E_0&=\frac{m}{\Sigma^2},\ \ J_3=\frac{m \kappa a}{\Sigma^2},\  \ q=\frac{e}{\Sigma},\ \ b_3=\frac{4\kappa a e}{3\Sigma},\\
J_1&=J_2=b_0=b_1=b_2=c_i=c'_i=0,\ \ i=1,2,3.
 \end{aligned}
 \eeQ

\mysection{Appendix}
The Killing vectors of $AdS$ spacetime are\\
\beQ \begin{aligned}
U_{40}=&\kappa^{-1} \frac{\partial}{\partial t},\\
U_{10}=&\kappa^{-1}\sin\theta \cos\psi\frac{\partial}{\partial r}
+\coth(\kappa r)\Big(\cos\theta\cos\psi\frac{\partial}{\partial \theta}-\frac{\sin\psi}{\sin \theta}\frac{\partial}{\partial \psi}\Big),\\
U_{20}=&\kappa^{-1}\sin\theta \sin\psi\frac{\partial}{\partial r} +\coth(\kappa r)\Big(\cos\theta\sin\psi\frac{\partial}{\partial
\theta}+\frac{\cos\psi}{\sin \theta}\frac{\partial}{\partial \psi}\Big),\\
U_{30}=&\kappa^{-1}\cos\theta\frac{\partial}{\partial r} -\coth(\kappa
r)\sin\theta\frac{\partial}{\partial \theta},\\
U_{14}=&\kappa^{-1}\tanh(\kappa r) \sin\theta \cos\psi\ \frac{\partial}{\partial t},\\
U_{24}=&\kappa^{-1}\tanh(\kappa r)\sin\theta \sin\psi \frac{\partial}{\partial t},\\
U_{34}=&\kappa^{-1}\tanh(\kappa r)\cos\theta
\frac{\partial}{\partial t},\\
U_{23}=&-\sin\psi\frac{\partial}{\partial \theta}-\frac{\cos\theta \cos \psi}{\sin \theta}\frac{\partial}{\partial
\psi},\\
U_{31}=&\cos\psi\frac{\partial}{\partial \theta}-\frac{\cos\theta \sin\psi}{\sin \theta}\frac{\partial}{\partial \psi},\\
U_{12}=&\frac{\partial}{\partial \psi}.
\end{aligned} \eeQ
We set $V_1=U_{23}$, $V_2=U_{31}$ and $V_3=U_{12}$ for convenience.

\ac
The authors would like to thank Professor X.Zhang for his suggestions and helpful discussions.
This work is done during the visit of the second author to Academy of Mathematics and System Science, Chinese Academy of Sciences.
He would like to thank the academy for its hospitality.

\end{document}